\documentclass[12pt]{amsart}
\usepackage{epsfig,amscd,amssymb,times}
\newtheorem{theorem}{Theorem}
\newtheorem{lemma}[theorem]{Lemma}
\newtheorem{proposition}[theorem]{Proposition}
\newtheorem{corollary}[theorem]{Corollary}
\theoremstyle{definition}
\newtheorem{definition}[theorem]{Definition}
\newtheorem{problem}[theorem]{Problem}

\theoremstyle{remark}
\newtheorem{remark}[theorem]{Remark}

\def\varph{{\varphi}}
\def\ra{{\rightarrow}}
\def\lra{{\longrightarrow}}
\def\om{{\omega}}

\def\la{{\lambda}}
\def\al{{\alpha}}

\def\La{{\Lambda}}
\def\ga{{\gamma}}

\def\bZ{{\mathbb Z}}
\def\bQ{{\mathbb Q}}
\def\bR{{\mathbb R}}
\def\inv{{^{-1}}}
\def\Sg{{\Sigma_g}}
\def\Mg{{\mathcal M}_g}
\def\Msg{{\mathcal M}_{g,*}}
\def\M1g{{\mathcal M}_{g,1}}
\def\I1g{{\mathcal I}_{g,1}}

\newcommand\Hom{\operatorname{Hom}}

\newcommand\sgn{\operatorname{sgn}}
\newcommand\Flux{\operatorname{Flux}}
\newcommand\Symp{\operatorname{Symp}}
\newcommand\Ham{\operatorname{Ham}}

\newcommand\Cal{\operatorname{Cal}}
\newcommand\Ker{\operatorname{Ker}}

\newcommand\im{\operatorname{Im}}
\newcommand\Int{\operatorname{Int}}

\newcommand\BS{\operatorname{BSymp}}
\newcommand\ES{\operatorname{ESymp}}

\newcommand\BG{\operatorname{B\bar\Gamma}}
\newcommand\Diff{\operatorname{Diff}}
\newcommand\id{\operatorname{id}}
\catcode`\@=\active 

\begin{document}
\title[Characteristic classes of foliated surface bundles]
{Characteristic classes of foliated surface bundles with area-preserving 
holonomy}
\author{D.~Kotschick}
\address{Mathematisches Institut, Ludwig-Maxi\-mi\-lians-Universit\"at 
M\"unchen,
Theresienstr.~39, 80333 M\"unchen, Germany}
\email{dieter{\char'100}member.ams.org}
\author{S.~Morita}
\address{Department of Mathematical Sciences\\
University of Tokyo \\Komaba, Tokyo 153-8914\\
Japan}
\email{morita{\char'100}ms.u-tokyo.ac.jp}

\thanks{The first author is grateful to the {\sl Deutsche 
Forschungsgemeinschaft}
for support of this work. The second author is partially supported by JSPS 
Grant
16204005}

\begin{abstract}
Making use of the {\it extended} flux homomorphism defined in~\cite{KM03}
on the group $\Symp\Sg$ of symplectomorphisms of a closed oriented surface
$\Sg$ of genus $g\geq 2$, we introduce new characteristic classes of foliated
surface bundles with symplectic, equivalently area-preserving, total holonomy.
These characteristic classes are stable with respect to $g$ and we show that
they are highly non-trivial. We also prove that the second homology of
the group $\Ham\Sg$ of Hamiltonian symplectomorphisms of $\Sg$, equipped with
the discrete topology, is very large for all $g\geq 2$.
\end{abstract}

\maketitle

\section{Introduction}\label{s:intro}

In this paper we study the homology of symplectomorphism groups of
surfaces considered as discrete groups. We shall prove that certain
homology groups are highly non-trivial by constructing characteristic
classes of foliated surface bundles with area-preserving holonomy, and
proving non-vanishing results for them.

Let $\Sg$ be a closed oriented surface of genus $g\geq 2$, and $\Diff_+\Sg$
its group of orientation preserving selfdiffeomorphisms. We fix an area
form $\omega$ on $\Sg$, which, for dimension reasons, we can also think of
as a symplectic form. We denote by $\Symp\Sg$ the subgroup of $\Diff_+\Sg$
preserving the form $\om$. The classifying space $\BS ^\delta\Sg$ for the group
$\Symp\Sg$ with the {\it discrete} topology is an Eilenberg-MacLane space
$K(\Symp^\delta\Sg,1)$ which classifies foliated $\Sg$-bundles with
area-preserving total holonomy groups.

Our construction of characteristic classes proceeds as follows.
Let $\Symp_0\Sg$ be the identity component of $\Symp\Sg$.
A well-known theorem of Moser~\cite{Moser} concerning volume-preserving
diffeomorphisms implies that the quotient $\Symp\Sg/\Symp_0\Sg$ can be
naturally identified with the mapping class group $\Mg$, so that
we have an extension
$$
1\lra\Symp_0\Sg\lra\Symp\Sg\overset{p}{\lra}\Mg\lra 1 \ .
$$
There is a surjective homomorphism $\Flux\colon\Symp_0\Sg\ra H^1(\Sg;\bR)$,
called the flux homomorphism. In~\cite{KM03} we proved that this homomorphism
can be extended to a crossed homomorphism
$$
\widetilde\Flux\colon\Symp\Sg\lra H^1(\Sg;\bR) \ ,
$$
which we call the {\it extended} flux homomorphism. This extension is 
essentially
unique in the sense that the associated cohomology class
$$
[\widetilde\Flux]\in H^1(\Symp^\delta\Sg;H^1(\Sg;\bR))
$$
with twisted coefficients is uniquely defined. Now we consider the powers
$$
[\widetilde\Flux]^k\in H^k(\Symp^\delta\Sg;H^1(\Sg;\bR)^{\otimes k})\quad
(k=2,3,\cdots) \ ,
$$
and apply $\Mg$-invariant homomorphisms
$$
\la\colon H^1(\Sg;\bR)^{\otimes k}\lra \bR
$$
to obtain cohomology classes
$$
\la([\widetilde\Flux]^k)\in H^{k}(\BS^\delta\Sg;\bR)
$$
with constant coefficients.
The usual cup-product pairing $H^1(\Sg;\bR)^{\otimes 2}\ra\bR$
is the main example of such a homomorphism $\la$.

This method of constructing constant cohomology classes out of twisted ones
was already used in~\cite{Morita96} in the case of the mapping class group,
where the Torelli group (respectively the Johnson homomorphism) played the
role of $\Symp_0\Sg$ (respectively of the flux homomorphism) here. In that
case, it was proved in {\it loc.~cit.}~that all the Mumford--Morita--Miller
classes can be obtained in this way. The precise formulae were given
in~\cite{KM96,KM01}, with the important implication that no other classes
appear. In our context here, we can go further by enhancing the coefficients
$\bR$ to associated $\bQ$-vector spaces which appear as the targets of various
multiples of the {\it discontinuous} cup-product pairing
$$
H^1(\Sg;\bR)\otimes_\bZ H^1(\Sg;\bR)\lra S^2_\bQ\bR \ ,
$$
where $S^2_\bQ\bR$ denotes the second symmetric power of $\bR$ over
$\bQ$, see Section~\ref{s:main} for the details. In this way, we obtain
many new characteristic classes in
$$
H^*(\Symp^\delta\Sg;S^*(S^2_\bQ\bR)) \ ,
$$
where
$$
S^*(S^2_\bQ\bR)=\bigoplus_{k=1}^\infty S^k(S^2_\bQ\bR)
$$
denotes the symmetric algebra of $S^2_\bQ\bR$.
On the other hand, we proved in~\cite{KM03} that any power
$$
e_1^k\in H^{2k}(\Symp^\delta\Sg;\bQ)
$$
of the first Mumford--Morita--Miller class $e_1$ is non-trivial for
$k\leq \frac{g}{3}$. Now we can consider the cup products of $e_1^k$ with
the new characteristic classes defined above. The main purpose of the present
paper is to prove that these characteristic classes are all non-trivial in a
suitable stable range.

The contents of this paper is as follows. In Section~\ref{s:main} precise
statements of the main results are given. In Section~\ref{s:trans} we study,
in detail, the transverse symplectic class of foliated $\Sg$-bundles with
area-preserving total holonomy groups. In Section~\ref{s:open} we construct
two kinds of {\it extended} flux homomorphisms for open surfaces
$\Sigma_g^0=\Sg\setminus D^2$. We compare these extended flux homomorphisms
with the one obtained in the case of closed surfaces. This is used in
Section~\ref{s:secondopen} to generalize our result on the second homology of
the symplectomorphism group to the case of open surfaces. Sections~\ref{s:open}
and~\ref{s:secondopen} are the heart of this paper. Then in Section~\ref{s:proofs} we use the results of the previous sections to show the non-t 
riviality of cup
products of various characteristic classes, thus yielding proofs of the main
results about the homology of symplectomorphism groups as discrete groups. In
the final Section~\ref{s:misc} we give definitions of yet more characteristic
classes, other than the ones given in Section~\ref{s:main}. We propose several
conjectures and problems about them.

\section{Statement of the main results}\label{s:main}

Consider the usual cup-product pairing
$$
\iota\colon H^1(\Sg;\bR)\otimes H^1(\Sg;\bR)\lra \bR
$$
in cohomology, dual to the intersection pairing in homology. For simplicity,
we denote $\iota(u,v)$ by $u\cdot v$, where $u,v\in H^1(\Sg;\bR)$.
We first lift this pairing as follows.

As in Section~\ref{s:intro}, let $S^2_{\bQ}\bR$ denote the second symmetric
power of $\bR$ over $\bQ$. In other words, this is a vector space over $\bQ$,
consisting of the homogeneous polynomials of degree two generated by the
elements of $\bR$ considered as a vector space over $\bQ$. For each element
$a\in\bR$, we denote by $\hat a$ the corresponding element in $S^1_\bQ\bR$.
Thus any element in $S^2_\bQ\bR$ can be expressed as a finite sum
$$
\hat a_1\hat b_1+\cdots +\hat a_k\hat b_k
$$
with $a_i, b_i\in\bR$. We have a natural projection
$$
S^2_\bQ\bR\lra \bR
$$
given by the correspondence $\hat a\mapsto a\ (a\in\bR)$.

With this terminology we make the following definition:
\begin{definition}[{\it Discontinuous intersection pairing}]
Define a pairing
$$
\tilde\iota\colon H^1(\Sg;\bR)\times H^1(\Sg;\bR)\lra S^2_\bQ\bR
$$
as follows. Choose a basis $x_1,\cdots , x_{2g}$ of $H^1(\Sg;\bQ)$.
For any two elements, $u,v\in H^1(\Sg;\bR)$, write
$$
u=\sum_i a_i x_i,\quad v=\sum_i b_i x_i\quad (a_i, b_i\in\bR) \ .
$$
Then we set
$$
\tilde\iota (u,v)=\sum_{i,j} \iota(x_i,x_j) \hat a_i \hat b_j \in S^2_\bQ
\bR \ .
$$
\end{definition}
Clearly $\tilde\iota$ followed by the projection $S^2_\bQ\bR\lra \bR$ is
nothing but the usual intersection or cup-product pairing $\iota$.
Henceforth we simply write $u\odot v$ for $\tilde\iota (u,v)$.

It is easy to see that $\tilde\iota$ is well defined independently of the
choice of basis in $H^1(\Sg;\bQ)$. We can see from the following proposition
that $\tilde\iota$ enumerates all $\bZ$-multilinear skew-symmetric pairings
on $H^1(\Sg;\bR)$ which are $\Mg$-invariant. Let $\La^2_{\bZ} H^1(\Sg;\bR)$
denote the second exterior power, over $\bZ$, of $H^1(\Sg;\bR)$ considered as
an abelian group, rather than as a vector space over $\bR$. Also let
$\left(\La^2_{\bZ} H^1(\Sg;\bR)\right)_{\Mg}$ denote the abelian group of
coinvariants of $\La^2_{\bZ} H^1(\Sg;\bR)$ with respect to the natural action
of $\Mg$.
\begin{proposition}\label{prop:la2}
There exists a canonical isomorphism
$$
\left(\La^2_{\bZ} H^1(\Sg;\bR)\right)_{\Mg}\cong S^2_{\bQ} \bR
$$
given by the correspondence
$$
\left(\sum_i a_i u_i\right)\land
\left(\sum_j b_j v_j\right)\longmapsto
\sum_{i,j} \iota(u_i,v_j)\ \hat a_i \hat b_j\ \ ,
$$
where $a_i, b_j\in\bR$, $u_i,v_j\in H^1(\Sg;\bQ)$.
\end{proposition}
In order not to digress, we refer the reader to the appendix for a proof.

Before we can define some new cocycles on the group $\Symp\Sg$, we have to
recall some facts from~\cite{KM03}. The symplectomorphism group $\Symp\Sg$
acts on its identity component by conjugation, and acts on $H^{1}(\Sg;\bR)$
from the left via $\varph(w)=(\varph^{-1})^*(w)$.
The flux homomorphism $\Flux\colon\Symp_{0}\Sg\rightarrow H^1(\Sg;\bR)$
is equivariant with respect to these actions by Lemma~6 of~\cite{KM03}.
Its extension $\widetilde{\Flux}\colon\Symp\Sg\rightarrow H^1(\Sg;\bR)$
is a crossed homomorphism for the above action in the sense that
\begin{equation}\label{eq:eflux}
\widetilde{\Flux}(\varph\psi)=\widetilde{\Flux}(\varph)+
(\varph\inv)^*\widetilde{\Flux}(\psi) \ .
\end{equation}
\begin{definition}\label{def:al}
     Let $\varph_1,\ldots,\varph_{2k}\in\Symp\Sg$, and
     $$
\xi_i=((\varph_1 \ldots\varph_{i-1})^{-1})^*\widetilde{\Flux}(\varph_i) \ .
$$
Define a $2k$-cocycle $\tilde\al^{(k)}$ with values in $S^{k}(S^2_\bQ\bR)$ by
\begin{align*}
\tilde\al^{(k)}&(\varph_1,\ldots,\varph_{2k})\\
&=
\frac{1}{(2k)!}\sum_{\sigma\in\mathfrak{S}_{2k}}
\sgn\sigma\ (\xi_{\sigma(1)}\odot\xi_{\sigma(2)})\ldots
(\xi_{\sigma(2k-1)}\odot\xi_{\sigma(2k)})\in S^{k}(S^2_\bQ\bR) \ ,
\end{align*}
where the sum is over permutations in the symmetric group $\mathfrak{S}_{2k}$.
\end{definition}

That $\tilde\al^{(k)}$ is indeed a cocycle is easy to check by a standard
argument in the theory of cohomology of groups using~\eqref{eq:eflux}. Thus
we have the corresponding cohomology classes
$$
\tilde\al^{(k)}\in H^{2k}(\Symp^\delta\Sg;S^{k}(S^2_\bQ\bR)) \ ,
$$
denoted by the same letters. If we apply the canonical projection
$S^{k}(S^2_\bQ\bR)\ra\bR$ to these classes, we obtain real cohomology classes
$$
\al^k\in H^{2k}(\Symp^\delta\Sg;\bR) \ ,
$$
which are the usual cup products of the first one $\al\in 
H^2(\Symp^\delta\Sg;\bR)$.
The refined classes $\tilde\al^{(k)}$ can be considered as a twisted version of
{\it discontinuous invariants} in the sense of~\cite{Morita85} arising from the
flux homomorphism.

Now we can state our first main result.
\begin{theorem}\label{th:hev}
For any $k\geq 1$ and $g\geq 3k$, the characteristic classes
$$
e_1^k, e_1^{k-1}\tilde\al, \ldots, e_1\tilde\al^{(k-1)}, \tilde\al^{(k)}
$$
induce a surjective homomorphism
$$
H_{2k}(\Symp^{\delta}\Sg;\bZ)\lra
\bZ\oplus S^{2}_{\bQ}\bR\oplus\cdots\oplus S^{k}(S^2_\bQ\bR) \ .
$$
\end{theorem}
For $k=1$ this is not hard to see, so we give the proof right away.
For $k>1$ the proof is given in Section~\ref{s:proofs} below and requires
the technical results developed in the body of this paper.

Consider the subgroup $\Ham\Sg$ of $\Symp_0\Sg$ consisting of all
Hamiltonian symplectomorphisms of $\Sg$. As is well known
(see~\cite{Banyaga,MS}), we have an extension
\begin{equation}\label{Ham}
1\lra\Ham\Sg\lra\Symp_0\Sg\ \overset{\Flux}{\lra}\
H^1(\Sg;\bR)\lra 1 \ .
\end{equation}
This gives rise to a $5$-term exact sequence in cohomology:
\begin{align*}
0 \lra H^{1}(H^1(\Sg;\bR)^\delta;\bZ) &\stackrel{\Flux^{*}}{\lra}
H^{1}(\Symp_0^\delta\Sg;\bZ)\lra
H^{1}(\Ham^\delta\Sg;\bZ)^{H^1_{\bR}}\\
&\lra H^{2}(H^1(\Sg;\bR)^\delta;\bZ)\stackrel{\Flux^{*}}{\lra}
H^{2}(\Symp_0^\delta\Sg;\bZ) \ ,
\end{align*}
where we have written $H^1_{\bR}$ for $H^1(\Sg;\bR)$. Now $\Ham\Sg$ is
a perfect group by a result of Thurston~\cite{Thurston0}, see also
Banyaga~\cite{Banyaga}. Therefore, $\Flux^{*}$ injects the second
cohomology of $H^1(\Sg;\bR)$ as a discrete group into that of
$\Symp_0^\delta\Sg$. By definition, the class $\tilde\alpha$ is the image
of the class of $\tilde\iota$ under $\Flux^{*}$. So $\tilde\alpha$ is
nontrivial on $\Symp_{0}^{\delta}\Sg$, and is defined on the whole
$\Symp^{\delta}\Sg$. We conclude that $\tilde\alpha$ defines a
surjective homomorphism
$$
H_{2}(\Symp^{\delta}\Sg;\bZ)\rightarrow S^{2}_{\bQ}\bR \ ,
$$
for any $g\geq 2$. We already proved in~\cite{KM03} that the first
Mumford--Morita--Miller class $e_{1}$ defines a surjection
$H_{2}(\Symp^{\delta}\Sg;\bZ)\longrightarrow \bZ$ for all $g\geq 3$.
Clearly the two classes are linearly independent because $\tilde\alpha$
is nonzero on $\Symp_{0}^{\delta}\Sg$, to which $e_{1}$ restricts trivially.
This proves Theorem~\ref{th:hev} in the easy case when $k=1$.

We can restrict the homomorphism
$$
\Flux^*\colon H^*(H^1(\Sg;\bR)^\delta;\bR)\lra
H^*(\Symp_0^\delta\Sg;\bR)
$$
to the {\it continuous cohomology}
$$
H^*_{ct}(H^1(\Sg;\bR)^\delta;\bR)\cong H_*(T^{2g};\bR)
\subset H^*(H^1(\Sg;\bR)^\delta;\bR) \ ,
$$
see Section~\ref{s:trans} for the precise definition. Thereby we obtain a ring
homomorphism
$$
\Flux^*\colon H_*(T^{2g};\bR)\longrightarrow H^*(\Symp_0^\delta\Sg;\bR) \ ,
$$
where $T^{2g}=K(\pi_1\Sg,1)$ is the Jacobian manifold of $\Sg$, and the ring
structure on the homology of $T^{2g}$ is induced by the Pontrjagin product.
Let $\om_0\in H_2(T^{2g};\bR)$ be the homology class represented by the dual
of the standard symplectic form on $T^{2g}$. We decompose the
$Sp(2g,\bR)$-module $H_k(T^{2g};\bR)$ into irreducible components.
For this, consider the homomorphism
$$
\om_0\wedge \ \colon H_{k-2}(T^{2g};\bR)\longrightarrow H_k(T^{2g};\bR)
$$
induced by the wedge product with $\om_0$. On the one hand, it is easy
to see using Poincar\'e duality on $T^{2g}$, that the above homomorphism is
surjective for any $k\geq g+1$. On the other hand, it is well-known
(see~\cite{FH}), that the kernel of the contraction homomorphism
$$
C\colon H_k(T^{2g};\bQ)\longrightarrow H_{k-2}(T^{2g};\bQ)
$$
induced by the intersection pairing
$H_2(T^{2g};\bQ)\cong \La^2 H_1(\Sg;\bQ)\ra\bQ$
is the irreducible representation of the algebraic group $Sp(2g,\bQ)$
corresponding to the Young diagram $[1^{k}]$ for any $k\leq g$.
Let $[1^k]_\bR=[1^k]\otimes\bR$ denote the real form of this representation.
Then we have a direct sum decomposition
\begin{equation}\label{eq:irrep}
H_k(T^{2g};\bR)= [1^k]_\bR\oplus \om_0\land H_{k-2}(T^{2g};\bR)
\quad (k\leq g) \ .
\end{equation}
\begin{theorem}\label{th:om2}
The kernel of the homomorphism
$$
\Flux^*\colon H_*(T^{2g};\bR)\lra H^*(\Symp_0^\delta\Sg;\bR)
$$
induced by the flux homomorphism is the ideal generated by the subspace
$\om_0\land H_1(T^{2g};\bR)\subset H_3(T^{2g};\bR)$, and the image of this
homomorphism can be described as
$$
\im \Flux^*\cong \bR\oplus \bigoplus_{k=1}^{g} [1^k]_\bR \ ,
$$
where $\bR$ denotes the image of the subspace of $H_2(T^{2g};\bR)$ spanned
by $\om_0$.
\end{theorem}

The group $H^1(\Sg;\bR)$ acts on $\Ham\Sg$ by outer automorphisms.
In particular, it acts on the homology $H_*(\Ham^{\delta}\Sg;\bZ)$
of the discrete group $\Ham^{\delta}\Sg$ so that we can consider the
coinvariants $H_*(\Ham^{\delta}\Sg;\bZ)_{H^1_{\bR}}$, where for simplicity
we have written $H^1_{\bR}$ instead of $H^1(\Sg;\bR)$.
The following result shows that the second homology group
$H_2(\Ham^{\delta}\Sg;\bZ)$ is highly non-trivial.
\begin{theorem}\label{th:ham}
For any $g\geq 2$, there exists a natural injection
$$
H^1(\Sg;\bR)\subset H_2(\Ham^{\delta}\Sg;\bZ)_{H^1_{\bR}} \ .
$$
\end{theorem}

\section{The transverse symplectic class}\label{s:trans}

Since $\Sg$ is an Eilenberg-MacLane space, the total space
of the universal foliated $\Sg$-bundle over the classifying
space $\BS^\delta\Sg$ is again a $K(\pi,1)$ space.
Hence if we denote by $\ES^{\delta}\Sg$ the fundamental group of
this total space, then we obtain a short
exact sequence
\begin{equation}\label{eq:esymp}
1\lra\pi_1\Sg\lra\ES^{\delta}\Sg\lra\Symp^{\delta}\Sg\lra 1
\end{equation}
and any cohomology class of the total space can be considered as
an element in the group cohomology of $\ES^{\delta}\Sg$.
Now on the total space of any foliated $\Sg$-bundle with total holonomy
group contained in $\Symp\Sg$ there is a closed $2$-form $\tilde\om$
which restricts to the symplectic form $\om$ on each fiber. At the
universal space level, the de Rham cohomology class of $\tilde\om$
defines a class $v\in H^2(\ES^\delta;\bR)$ which we call the transverse
symplectic class. We normalize the symplectic form $\om$ on $\Sg$
so that its total area is equal to $2g-2$. It follows that the
restriction of $v$ to a fiber is the same as the negative of the Euler
class $e\in H^2(\ES^\delta;\bR)$ of the vertical tangent bundle.

Let $\ES_0\Sg$ denote the subgroup of $\ES^{\delta}\Sg$ obtained by
restricting the extension~\eqref{eq:esymp} to $\Symp_0\Sg\subset\Symp\Sg$.
Since any foliated $\Sg$-bundle with total holonomy in $\Symp_0\Sg$ is
trivial as a differentiable $\Sg$-bundle, there exists an isomorphism
$$
\ES_0\Sg\cong \pi_1\Sg\times\Symp_0\Sg \ .
$$
Henceforth we identify the above two groups.
By the K\"unneth decomposition, we have an isomorphism
\begin{equation}\label{eq:kunneth}
\begin{split}
H^2(&\ES_0^\delta\Sg;\bR)\cong H^2(\Sg;\bR)\oplus\\
&\oplus \left(H^1(\Sg;\bR)\otimes
H^1(\Symp_0^\delta\Sg;\bR)\right)\oplus H^2(\Symp_0^\delta\Sg;\bR) \ ,
\end{split}
\end{equation}
where we identify $H^*(\pi_1\Sg;\bR)$ with $H^*(\Sg;\bR)$.
Let $\mu\in H^2(\Sg;\bZ)$ be the fundamental cohomology class
of $\Sg$. Clearly the Euler class $e\in H^2(\ES_0^\delta\Sg;\bR)$
is equal to $(2-2g)\mu$.
The flux homomorphism 
gives rise to an element
\begin{align*}
[\Flux]\in \hspace{3mm} &\Hom_\bZ (H_1(\Symp_0^\delta\Sg;\bZ),H^1(\Sg;\bR))\\
\cong & \Hom_\bR(H_1(\Symp_0^\delta\Sg;\bR),H^1(\Sg;\bR))\\
\cong & H^1(\Sg;\bR)\otimes H^1(\Symp_0^\delta\Sg;\bR) \ ,
\end{align*}
where the last isomorphism exists because $H^1(\Sg;\bR)$
is finite dimensional. Choose a symplectic basis
$x_1,\ldots,x_g$, $y_1,\ldots,y_g$
of $H_1(\Sg;\bR)$ and denote by $x_1^*,\ldots,x_g^*$, $y_1^*,\ldots,y_g^*$
the dual basis of $H^1(\Sg;\bR)$. Then Poincar\'e duality
$H_1(\Sg;\bR)\cong H^1(\Sg;\bR)$ is given by the
correspondence $x_i\mapsto -y_i^*, y_i\mapsto x_i^*$.
The element $[\Flux]$ can be described explicitly as
\begin{equation}\label{eqn:flux}
[\Flux]=\sum_{i=1}^{g} (x_i^*\otimes \tilde x _i+y_i^*\otimes \tilde y_i)
\in H^1(\Sg;\bR)\otimes H^1(\Symp_0^\delta\Sg;\bR)
\end{equation}
where $\tilde x_i,\tilde y_i\in H^1(\Symp_0^\delta\Sg;\bR)\cong
\Hom (H_1(\Symp_0\Sg;\bZ);\bR)$ is defined by the
equality
$$
\Flux(\varph)=\sum_{i=1}^{g}
(\tilde x_i(\varph) x_i+\tilde y_i(\varph) y_i)\quad (\varph\in \Symp_0\Sg).
$$
The elements $\tilde x_i,\tilde y_i$ can be also interpreted
as follows.
The flux homomorphism induces a homomorphism in cohomology
\begin{equation}\label{eq:flux*}
\Flux^*\colon H^*(H^1(\Sg;\bR)^\delta;\bR)\lra H^*(\Symp_0^\delta\Sg;\bR)
\end{equation}
where the domain
$$
H^*(H^1(\Sg;\bR)^\delta;\bR)\cong \Hom_{\bZ}(\La_{\bZ}^*(H^1(\Sg;\bR)),\bR)
$$
is the cohomology group of $H^1(\Sg;\bR)$ considered as a {\it discrete
abelian group}, rather than as a vector space over $\bR$, so that it is
a very large group. Its {\it continuous part} is defined as
\begin{align*}
H^*_{ct}(H^1(\Sg;\bR)^\delta;\bR)=&\Hom_{\bR}(\La_{\bR}^*(H^1(\Sg;\bR)),\bR)\\
\subset &\Hom_{\bZ}(\La_{\bR}^*(H^1(\Sg;\bR)),\bR)\\
\subset &\Hom_{\bZ}(\La_{\bZ}^*(H^1(\Sg;\bR)),\bR)\cong
H^*(H^1(\Sg;\bR)^\delta;\bR) \ ,
\end{align*}
where the second inclusion is induced by the natural projection
$$
\La_{\bZ}^*(H^1(\Sg;\bR))\longrightarrow \La_{\bR}^*(H^1(\Sg;\bR)) \ .
$$
Denoting by $T^{2g}=K(\pi_1\Sg,1)$ the Jacobian torus of $\Sg$, there is
a canonical isomorphism
$$
\La^*_{\bR}H^1(\Sg;\bR)\cong H^*(T^{2g};\bR) \ ,
$$
so that we can identify
$$
H^*_{ct}(H^1(\Sg;\bR)^\delta;\bR)\cong \Hom_\bR(H^*(T^{2g};\bR);\bR)
\cong H_*(T^{2g};\bR) \ .
$$
Thus, by restricting the homomorphism $\Flux^{*}$ in~\eqref{eq:flux*}
to the continuous cohomology, we obtain a homomorphism
$$
\Flux^*\colon H_*(T^{2g};\bR)\longrightarrow H^*(\Symp_0^\delta\Sg;\bR) \ .
$$
It is easy to see that under this homomorphism we have
$$
\tilde x_i=\Flux^*(x_i) \ ,\quad \tilde y_i=\Flux^*(y_i) \ .
$$
Let $\om_0\in \La^2_{\bR}H_1(\Sg;\bR)$ be the symplectic class
defined by
$$
\om_0=\sum_{i=1}^g x_i\land y_i
$$
and set
$$
\tilde\om_0=\Flux^*(\om_0)
=\sum_{i=1}^{g} \tilde x_i \tilde y_i\in H^2(\Symp_0^\delta\Sg;\bR) \ .
$$

\begin{lemma}\label{lem:flux2}
We have the equality
$$
[\Flux]^2=-2\mu\otimes \tilde\om_0\in
H^2(\Sg;\bR)\otimes H^2(\Symp_0^\delta\Sg;\bR) \ .
$$
\end{lemma}
\begin{proof}
A direct calculation using the expression~\eqref{eqn:flux} yields
$$
[\Flux]^2=-\sum_{i=1}^{g} (x_i^* y_i^*\otimes\tilde x_i\tilde y_i
+y_i^* x_i^*\otimes \tilde y_i\tilde x_i) \ .
$$
Since $x_i^*y_i^*=-y_i^*x_i^*=\mu$, we obtain
$$
[\Flux]^2=-2\mu\otimes\sum_{i=1}^{g} \tilde x_i\tilde y_i
=-2\mu\otimes\tilde\om_0
$$
as required.
\end{proof}

Now we can completely determine the transverse symplectic class $v$ of
foliated $\Sg$-bundles whose total holonomy groups are contained
in the identity component $\Symp_0\Sg$ of $\Symp\Sg$ as follows.
\begin{proposition}\label{prop:v}
On the subgroup $\ES_0^\delta\Sg$ the transverse symplectic class
$v\in H^2(\ES_0^\delta\Sg;\bR)$ is given by
$$
v=(2g-2)\mu+[\Flux]+\frac{1}{2g-2}\tilde\om_0
$$
under the isomorphism \eqref{eq:kunneth}.
Furthermore, the homomorphism
$$
\tilde\om_0\otimes H^1_{ct}(\Symp_0^\delta\Sg;\bR)\lra
H^3(\Symp_0^\delta\Sg;\bR)
$$
induced by the cup product is trivial, where
$H^1_{ct}(\Symp_0^\delta\Sg;\bR)\cong H_1(\Sg;\bR)$
denotes the subgroup of $H^1(\Symp_0^\delta\Sg;\bR)$
generated by the continuous cohomology classes
$\tilde x_i, \tilde y_i$.
In particular, $\tilde\om_0^2=0$.
\end{proposition}
\begin{proof}
Since the restriction of $v$ to each fiber is equal to the negative of
that of the Euler class $e$ by our normalization, $v$ restricts to
$(2g-2)\mu$ on the fiber. This gives the first component of the formula.
The second component follows from Lemma~8 of~\cite{KM03}.
Thus we can write
$$
v=(2g-2)\mu+[\Flux]+ \ga \in H^*(\Sg;\bR)\otimes H^*(\Symp_0^\delta\Sg;\bR)
$$
for some $\ga\in H^2(\Symp_0^\delta\Sg;\bR)$.
Now observe that $v^2=0$ because $\tilde\om^2=0$.
Also, because we have restricted to $\Symp_{0}\Sg$, we have
$\mu^2=0$, $\mu [\Flux]=0$. Hence we obtain
$$
[\Flux]^2+\ga^2+ 2(2g-2)\mu \ga+2 [\Flux] \ga=0 \ .
$$
It follows that
\begin{align*}
[\Flux]^2+2(2g-2)\mu \ga&=0\\
[\Flux] \ga&=0\\
\ga^2&=0
\end{align*}
because these three elements belong to different summands in the
K\"unneth decomposition
of $H^*(\Sg;\bR)\otimes H^*(\Symp_0^\delta;\bR)$.
If we combine the first equality above and Lemma~\ref{lem:flux2},
then we can conclude that
\begin{equation}\label{eqn:ga}
\ga=\frac{1}{2g-2}\tilde\om_0 \ .
\end{equation}
This proves the first claim of the proposition.
If we substitute~\eqref{eqn:flux} and~\eqref{eqn:ga} in
the second equality above, then we see that
$\tilde\om_0\tilde x_i=\tilde\om_0\tilde y_i=0$ for any $i$,
whence the second claim.
Observe that the third equality $\ga^2=0$, which is equivalent to
$\tilde\om_0^2=0$ by the above, is a consequence of the second claim.
\end{proof}

Now we can calculate the restriction of the cocycle $\alpha$ defined
in Section~\ref{s:main} to the identity component $\Symp_{0}\Sg$.
\begin{proposition}\label{prop:alom}
Let $i\colon\Symp_0^\delta\Sg\ra\Symp^\delta\Sg$ be the inclusion.
Then 
$$
i^*\al = 2 \tilde\om_0\in H^2(\Symp_0^\delta\Sg;\bR)
$$
and
$$
i^*\al^2=0 \ .
$$
\end{proposition}
\begin{proof}
Let $\varph, \psi\in \Symp_0\Sg$ be any two elements. Then, by the
definition of $\al$, see Definition~\ref{def:al} and the subsequent
discussion, we have
$$
\al(\varph,\psi)=\iota(\Flux(\varph), \Flux(\psi)) \ .
$$
By the definition of the cohomology classes $\tilde x_i, \tilde y_i$, we have
$$
\Flux(\varph)=\sum_{i=1}^{g}
(\tilde x_i(\varph) x_i+\tilde y_i(\varph) y_i)\ ,\quad
\Flux(\psi)=\sum_{i=1}^{g}
(\tilde x_i(\psi) x_i+\tilde y_i(\psi) y_i)\ .
$$
Hence we obtain
$$
\al(\varph,\psi)=\sum_{i=1}^g \left\{\tilde x_i(\varph) \tilde y_i(\psi)
-\tilde y_i(\varph) \tilde x_i(\psi)\right\}\ .
$$
Using the Alexander--Whitney cup product, we also have
$$
\tilde\om_0(\varph,\psi)=\sum_{i=1}^g \tilde x_i(\varph) \tilde y_i(\psi)
=-\sum_{i=1}^g \tilde x_i(\psi) \tilde y_i(\varph) \ .
$$
Thus $\al=2 \tilde\om_0$ as required. The last statement follows from
Proposition~\ref{prop:v}.
\end{proof}
We see from this proposition that $i^*\al^2=0$ in
$H^4(\Symp_0^\delta\Sg;\bR)$, while we will prove the
non-triviality of $\al^2\in H^4(\Symp^\delta\Sg;\bR)$, which is a special
case of Theorem~\ref{th:hev}.
One could say that the non-triviality of $\al^2$ is realized
by an interaction of the two groups $\Symp_0\Sg$ and $\Mg$.

\section{Extended flux homomorphisms for open surfaces}\label{s:open}

We consider the open surface $\Sigma_g^0=\Sigma_g\setminus D$ obtained
from $\Sg$ by removing a closed embedded disk $D\subset\Sg$. Let
$j\colon\Sigma_g^0\ra \Sg$ be the inclusion. We denote by $\Symp^c\Sigma_g^0$
the symplectomorphism group of $(\Sigma_g^0,j^*\om)$ with {\it compact 
supports}.
Hence the group $\Symp^c\Sigma_g^0$ can be considered as a subgroup of 
$\Symp\Sg$
with inclusion $j\colon\Symp^c\Sigma_g^0\ra\Symp\Sg$. Let $\Symp_0^c\Sigma_g^0$
be the identity component of $\Symp^c\Sigma_g^0$. Clearly $\Symp_0^c\Sigma_g^0$
is a subgroup of $\Symp_0\Sg$.

Let $\M1g$ denote the mapping class group of $\Sg$ {\it relative to}
the embedded disk $D^2\subset\Sg$, equivalently, the mapping class
group of the compact surface $\overline{\Sigma_g^0}$ with boundary.
We have a natural homomorphism $p\colon\Symp^c\Sigma_g^0\ra\M1g$
which is easily seen to be surjective. Moreover, Moser's theorem~\cite{Moser}
adapted to the present case (see~\cite{Tsuboi} for a general statement),
implies that the kernel of this surjection is precisely the group
$\Symp_0^c\Sigma_g^0$. We summarize the situation in the following
diagram:
\begin{equation}
\begin{CD}
1@>>> \Symp^{c}_{0}\Sigma_g^{0} @>{i^c}>>
\Symp^c\Sigma_g^{0}@>{p}>> \M1g @>>> 1\\
@. @V{j_0}VV @V{j}VV @V{q}VV  @.\\
1@>>> \Symp_{0}\Sigma_g @>{i}>> \Symp\Sigma_g @>>> \Mg @>>> 1
\end{CD}
\label{eqn:GG}
\end{equation}
where $q\colon\M1g\ra\Mg$ denotes the natural projection.

The restriction of the flux homomorphism
\begin{equation}\label{eqn:Flux}
\Flux\colon\Symp_0\Sg\lra H^1(\Sg;\bR)
\end{equation}
to the subgroup $\Symp_0^c\Sigma_g^0$, denoted $j^*\Flux$, can be described
as follows. The restriction $j^{*}\omega$ of the area form $\omega$
to the open surface $\Sigma_{g}^{0}$ is exact. Choose a $1$-form $\la$ 
$\Sigma_g^0$
such that $d\la=-j^*\om$. Then, for any element $\varph\in\Symp_0^c\Sigma_g^0$,
the $1$-form $\la-\varph^*\la$ is a closed form with compact support.
Hence the corresponding de Rham cohomology class $[\la-\varph^*\la]$,
which can be shown to be independent of the choice of $\la$,
is an element of the first cohomology group $H^1_c(\Sigma_g^0;\bR)$
of $\Sigma_g^0$ with compact support. It is easy to see that
$H^1_c(\Sigma_g^0;\bR)$ is canonically isomorphic to $H^1(\Sg;\bR)$,
and that under this isomorphism
\begin{equation}\label{eq:Fluxc}
(j^*\Flux)(\varph)=[\la-\varph^*\la]\in H^1_c(\Sigma_g^0;\bR)
\cong H^1(\Sg;\bR) 
\end{equation}
for all $\varph\in\Symp_0^c\Sigma_g^0$, see Lemma~10.14 of~\cite{MS}.
We obtain the following commutative diagram:
\begin{equation}
\begin{CD}
\Symp^{c}_{0}\Sigma_g^{0} @>{j^*\Flux}>> H^1_{c}(\Sigma_g^{0};\bR)\\
@VVV @VV{\cong}V\\
\Symp_{0}\Sigma_g @>>{\Flux}> H^1(\Sg;\bR).
\end{CD}
\label{eqn:FFlux}
\end{equation}
From now on we identify $H^1_{c}(\Sigma_g^{0};\bR)$ with $H^1(\Sg;\bR)$.

As was already mentioned in the Introduction, we proved in~\cite{KM03} that
the flux homomorphism~\eqref{eqn:Flux} can be extended to a crossed 
homomorphism
\begin{equation}\label{eqn:TFlux}
\widetilde{\Flux}\colon\Symp\Sg\lra H^1(\Sg;\bR) \ ,
\end{equation}
and that the extension is unique up to the addition of coboundaries.
The restriction of such a crossed homomorphism to the subgroup
$\Symp^c\Sigma_g^0\subset \Symp\Sg$, denoted $j^*\widetilde{\Flux}$,
is of course an extension of the flux homomorphism $j^*\Flux$.
However, we also have another extension of the same flux homomorphism
$j^*\Flux$ to the group $\Symp^c\Sigma_g^0$ as follows.
\begin{proposition}\label{prop:Fluxb}
The map
$$
\widetilde{\Flux}_c\colon\Symp^c\Sigma_g^0\lra H^1_c(\Sigma_g^0;\bR)
\cong H^1(\Sg;\bR) 
$$
defined by
$\widetilde{\Flux}_c(\varph)=[(\varph\inv)^*\la-\la]\in H^1_c(\Sigma_g^0;\bR)$
is a crossed homomorphism which extends the flux homomorphism $j^*\Flux$.
Its cohomology class
$[\widetilde{\Flux}_c]\in H^1(\Symp^c\Sigma_g^0;H^1(\Sg;\bR))$
is uniquely determined independently of the choice of the
$1$-form $\la$ such that $d\la=-j^*\om$.
\end{proposition}
\begin{proof}
First observe that, for any $\varph\in \Symp_0^c\Sigma_g^0$, we have
$$
[\la-\varph^*\la]=(\varph\inv)^*[\la-\varph^*\la]=[(\varph\inv)^*\la-\la]
$$
because $\varph$ acts trivially on $H^1(\Sg;\bR)$. Hence
by~\eqref{eq:Fluxc} we have
$$
\widetilde{\Flux}_c(\varph)=(j^*\Flux)(\varph) \ .
$$
Next, for any two elements $\varph, \psi\in \Symp^c\Sigma_g^0$, we have
\begin{align*}
\widetilde{\Flux}_c(\varph\psi)&=[((\varph\psi)\inv)^*\la-\la]\\
&=[(\varph\inv)^*\la-\la +(\varph\inv)^*(\psi\inv)^*\la-(\varph\inv)^*\la]\\
&=[(\varph\inv)^*\la-\la]+[(\varph\inv)^*((\psi\inv)^*\la-\la)]\\
&=\widetilde{\Flux}_c(\varph)+(\varph\inv)^*\widetilde{\Flux}_c(\psi).
\end{align*}
Therefore $\widetilde{\Flux}_c$ is a crossed homomorphism which extends 
$j^*\Flux$.

Finally, let $\la'$ be another $1$-form on $\Sigma_g^0$
satisfying $d\la'=-j^*\om$, and let $\widetilde{\Flux}'_c$ be the corresponding
crossed homomorphism. Then $a=\la'-\la$ is a closed $1$-form defining a de
Rham cohomology class $[a]\in H^1(\Sigma_g^0;\bR)\cong H^1(\Sg;\bR)$.
Now
\begin{align*}
\widetilde{\Flux}'_c(\varph)&=[(\varph\inv)^*\la'-\la']\\
&=[(\varph\inv)^*(\la+a)-(\la+a)]\\
&=[(\varph\inv)^*\la-\la]+[(\varph\inv)^*a-a]\\
&=\widetilde{\Flux}_c(\varph)+(\varph\inv)^*[a]-[a] \in H^1(\Sg;\bR)
\ .
\end{align*}
This shows that the difference $\widetilde{\Flux}'_c-\widetilde{\Flux}_c$
is a coboundary, completing the proof of the proposition.
\end{proof}

We have proved that the restriction $j^*\Flux$ of the flux
homomorphism~\eqref{eqn:Flux} to the subgroup $\Symp_0^c\Sigma_g^0$ has
two extensions to the group $\Symp^c\Sigma_g^0$ as a crossed homomorphism.
One is the restriction $j^*\widetilde{\Flux}$ of $\widetilde{\Flux}$, and
the other is $\widetilde{\Flux}_c$. We will show that these two crossed
homomorphisms are essentially different. More precisely, we will show that
the difference of these two crossed homomorphisms can be expressed by an
element of the cohomology group
\begin{equation}\label{eqn:hm1g}
H^1(\M1g;H^1(\Sg;\bR)) \ .
\end{equation}
It was proved in~\cite{Morita89} that $H^1(\M1g;H^1(\Sg;\bZ))$ is
isomorphic to $\bZ$ for all $g\geq 2$. It follows that the above
group~\eqref{eqn:hm1g} is isomorphic to $\bR$ and if
$$
k\colon \M1g \lra  H^1(\Sg;\bZ)
$$
is any crossed homomorphism whose cohomology class is a
generator of $H^1(\M1g;H^1(\Sg;\bZ))$, then the associated
crossed homomorphism
$$
k_\bR\colon\M1g \lra  H^1(\Sg;\bR)
$$
represents the element $1\in H^1(\M1g;H^1(\Sg;\bR))\cong\bR$.
Let $p^*k_\bR\in H^1(\Symp^c\Sg;H^1(\Sg;\bR))$ be the
class induced from $k_\bR$ by the projection $p\colon\Symp^c\Sg\ra\M1g$.

\begin{theorem}\label{th:h1sc}
For any $g\geq 2$, there exists an isomorphism
$$
H^1(\Symp^c\Sigma_g^0;H^1(\Sg;\bR))\cong \bR\oplus\Hom_{\bZ}(\bR,\bR)
$$
such that $p^*k_\bR$ and $[j^*\widetilde{\Flux}]$ represent
the classes $1\in \bR$ and $\id\in \Hom_{\bZ}(\bR,\bR)$ respectively.
\end{theorem}
\begin{proof}
The top extension in~\eqref{eqn:GG} gives rise to an exact sequence
\begin{align*}
& 0\ \lra \ H^1(\M1g;H^1(\Sg;\bR))\ \cong\bR\ \overset{p^*}{\lra} \\
   H^1(&\Symp^c\Sg;H^1(\Sg;\bR))\lra
H^1(\Symp^c_0\Sg;H^1(\Sg;\bR))^{\M1g}\lra\cdots .
\end{align*}
It was proved in~\cite{KM03} that $j^*\Flux$ induces an isomorphism
$$
j^*\Flux\colon H_1(\Symp^c_0\Sigma_g^0;\bZ)\cong
H^1_c(\Sigma_g^0;\bR)=H^1(\Sg;\bR) \ .
$$
Hence, by an argument similar to  the one given for $\Symp_0\Sg$ 
in~\cite{KM03},
we have an isomorphism
$$
H^1(\Symp^c_0\Sigma_g^0;H^1(\Sg;\bR))^{\M1g}\cong \Hom_\bZ(\bR,\bR) \ ,
$$
and clearly $j^*[\widetilde{\Flux}]$ corresponds to $\id\in \Hom_\bZ(\bR,\bR)$.
The result follows from this.
\end{proof}

\begin{theorem}\label{th:ffk}
We have the identity
$$
[\widetilde{\Flux}_c]=j^*[\widetilde{\Flux}]- p^*k_\bR
$$
in $H^1(\Symp^c\Sigma_g^0;H^1(\Sg;\bR))$.
\end{theorem}
\begin{proof}
Since both crossed homomorphisms $\widetilde{\Flux}_c$ and 
$j^*\widetilde{\Flux}$
are extensions of the flux homomorphism $j^*\Flux$, the proof of
Theorem~\ref{th:h1sc} implies that
$$
[\widetilde{\Flux}_c]=j^*[\widetilde{\Flux}]+a\ p^*k_\bR
$$
for some constant $a\in\bR$. Let $\I1g\subset\M1g$ denote the Torelli
subgroup consisting of mapping classes which act trivially on homology.
We set ${\mathcal I}\Symp^c\Sigma_g^0=p^{-1}(\I1g)\subset\Symp^c\Sigma_g^0$.
If we restrict the crossed homomorphims $\widetilde{\Flux}_c$ and
$\widetilde{\Flux}$ to this subgroup ${\mathcal I}\Symp^c\Sigma_g^0$, then
they become homomorphisms which depend only on the cohomology classes
$[\widetilde{\Flux}_c]$ and $j^*[\widetilde{\Flux}]$, and not on the
particular crossed homomorphisms representing these cohomology classes.
This is because any crossed homomorphism which is a coboundary is trivial
on ${\mathcal I}\Symp^c\Sigma_g^0$. It was proved in~\cite{Morita89}
that a generator of the group $H^1(\M1g;H^1(\Sg;\bZ))\cong \bZ$ is
characterized by the fact that the {\it Poincar\'e dual} of its value on
a single non-trivial element $\varph \in \I1g$ is equal to $\pm C 
\tau(\varph)$,
where $\tau\colon\I1g\ra \La^3 H_1(\Sg;\bZ)$ denotes the
(first) Johnson homomorphism and $C\colon\La^3 H_1(\Sg;\bZ)\ra H_1(\Sg;\bZ)$
denotes the contraction. Hence we only have to compute
the values of $\widetilde{\Flux}_c$ and $\widetilde{\Flux}$
on some particular element $\tilde\varph\in {\mathcal I}\Symp^c\Sigma_g^0$.
We choose such element as follows.

\begin{figure}
\epsfig{file=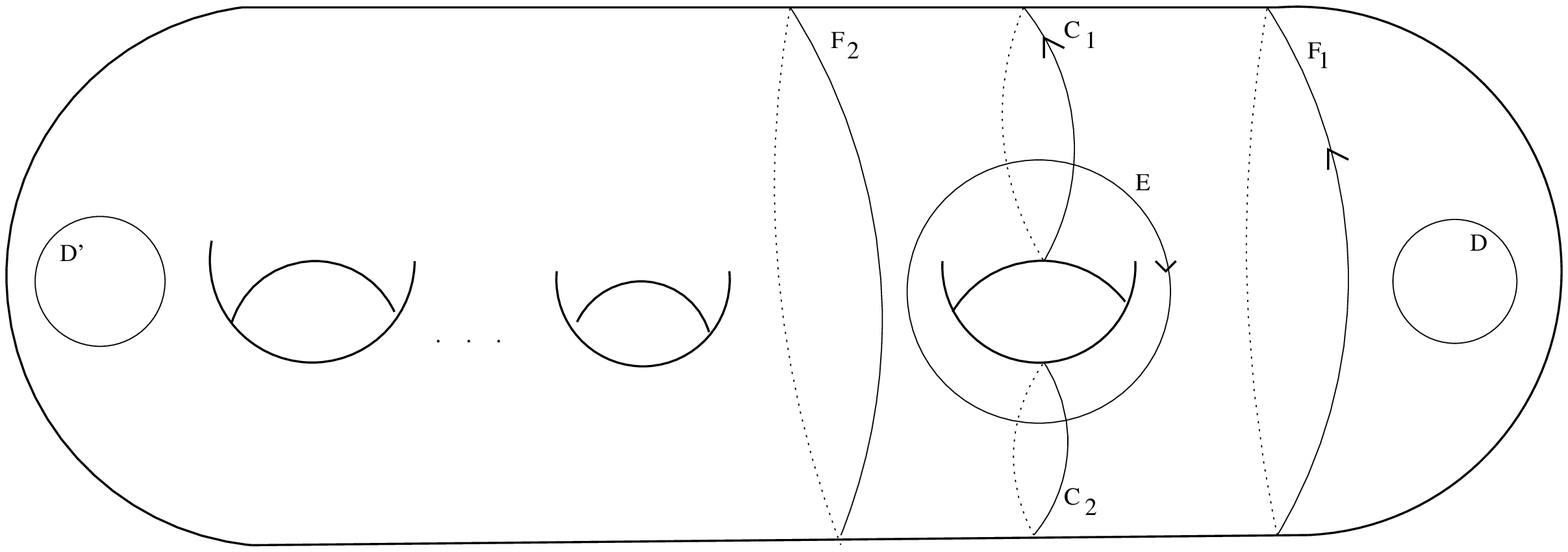,width=12cm}
\caption{ }\label{fig:torelli}
\end{figure}

Consider an embedded disk $D'$ in $\Sigma_g^0=\Sg\setminus D$ as depicted
in Figure~\ref{fig:torelli} and set $\Sigma_g^{00}=\Sg\setminus (D\cup D')$.
We also consider two disjoint simple closed curves $C_1, C_2\subset 
\Sigma_g^{00}$
as shown in Figure~\ref{fig:torelli}.
Let $\tau_i\ (i=1,2)$ denote the Dehn twist along $C_i$
and choose lifts $\tilde\tau_i\in\Symp^c\Sigma_g^0$ of $\tau_i$ with
supports in small neighborhoods of the $C_{i}$.
Now we set
$$
\tilde\varph=\tilde\tau_1 \tilde\tau_2^{-1} \ .
$$
Since $\varph=\tau_1\tau_2^{-1}$ belongs to the Torelli group $\I1g\subset 
\M1g$,
its lift $\tilde\varph$ is an element of the subgroup
${\mathcal I}\Symp^c\Sigma_g^0$. We would like to compute the difference
\begin{equation}\label{eq:difference}
\begin{split}
\widetilde{\Flux}_c(\tilde\varph)-\widetilde{\Flux}(\tilde\varph)
&= a\ p^*k_\bR(\tilde\varph)\\
&=a\ k(\varph)\in H^1_c(\Sigma_g^0;\bR)\cong H^1(\Sg;\bR) \ .
\end{split}
\end{equation}

Now we consider the long exact sequence
\begin{equation}\label{eq:s0}
\begin{split}
0\lra H^1(\Sg;\bR)&\overset{i^*}{\lra}H^1(\Sigma_g^{00};\bR)
\lra \\
& H^2(\Sg,\Sigma_g^{00};\bR)\cong\bR^2\lra
H^2(\Sg;\bR)\lra 0
\end{split}
\end{equation}
of the pair $(\Sg, \Sigma_g^{00})$.
Also consider the following short exact sequences
\begin{equation}\label{eq:s0'}
\begin{split}
0\lra H^1_c(\Sigma_g^0;\bR)\cong H^1(\Sg;\bR)&\overset{j^*}{\lra}
H^1(\Sigma_g^{00};\bR)\lra \bR\lra 0\\
0\lra H^1_c(\Sigma_g\setminus D';\bR)\cong H^1(\Sg;\bR)&\overset{(j')^*}{\lra}
H^1(\Sigma_g^{00};\bR)\lra \bR\lra 0 \ .
\end{split}
\end{equation}
In view of the above exact sequences~\eqref{eq:s0} and~\eqref{eq:s0'},
we can determine the value of~\eqref{eq:difference} in the
group $H^1(\Sigma_g^{00};\bR)$ because both groups
$H^1_c(\Sigma_g^0;\bR)$ and $H^1(\Sg;\bR)$ are embedded in it.
Choose $1$-forms
$\la$ and $\la'$ on $\Sigma_g^0=\Sg\setminus D$ and
$\Sg\setminus D'$ respectively, such that
$$
d\la=-j^*\om,\quad d\la'=-(j')^*\om
$$
where $j\colon\Sigma_g^0\ra\Sg$ and $j'\colon\Sg\setminus D'\ra\Sg$ are
the inclusions. Then $d(\la-\la')=0$ so that $\nu=\la-\la'$ is a
closed $1$-form on $\Sigma_g^{00}$. We have the identity
\begin{align*}
(\tilde\varph\inv)^*\la-\la&=(\tilde\varph\inv)^*(\la'+\nu)-(\la'+\nu)\\
&=(\tilde\varph\inv)^*\la'-\la'+(\tilde\varph\inv)^*\nu-\nu \ .
\end{align*}
Hence we have the equality
\begin{equation}\label{eq:[]}
[(\tilde\varph\inv)^*\la-\la]=[(\tilde\varph\inv)^*\la'-\la']+
[(\tilde\varph\inv)^*\nu-\nu]
\end{equation}
in the group $H^1(\Sigma_g^{00};\bR)$.
By the definition of $\widetilde{\Flux}_c$, we have
\begin{equation}\label{eq:fla}
\widetilde\Flux_c(\tilde\varph)=[(\tilde\varph\inv)^*\la-\la]\in
H^1_c(\Sigma_g^0;\bR)\subset H^1(\Sigma_g^{00};\bR) \ .
\end{equation}
Next we compute $\widetilde{\Flux}(\tilde\varph)$. For this, observe
that $\tilde\varph$ is isotopic to the identity as an element of
$\Symp^c(\Sg\setminus D')\subset \Symp\Sg$, although it is {\it not}
isotopic to the identity as an element of $\Symp^c(\Sg\setminus D)$.
Hence $\tilde\varph\in \Symp^c_0(\Sg\setminus D')$ and
$\widetilde{\Flux}(\tilde\varph)=((j')^*\Flux)(\tilde\varph)$.
By replacing $D$ with $D'$ in equation~\eqref{eq:Fluxc}, we obtain
\begin{equation}\label{eq:Fluxc'}
\begin{split}
((j')^*\Flux) &(\tilde\varph)=[\la'-\tilde\varph^*\la']\\
= &[(\tilde\varph\inv)^*\la'-\la']\in
H^1_c(\Sigma_g\setminus D';\bR)\subset H^1(\Sigma_g^{00};\bR) \ .
\end{split}
\end{equation}
By combining the equations~\eqref{eq:difference}, \eqref{eq:[]},
\eqref{eq:fla} and~\eqref{eq:Fluxc'}, we shall prove
\begin{equation}\label{eq:nutau}
[(\tilde\varph\inv)^*\nu-\nu]=
a k(\tilde\varph)=a k(\varph)=a PD\circ C \tau(\varph) \ ,
\end{equation}
where $PD$ denotes the Poincar\'e duality isomorphism.

Clearly $(\tilde\varph\inv)^*\nu-\nu$ is a closed $1$-form on $\Sg$ whose
support is contained in a neighborhood of $C_1\cup C_2$.
Hence the Poincar\'e dual of the de Rham cohomology class
$[(\tilde\varph\inv)^*\nu-\nu]\in H^1(\Sg;\bR)$ is a
multiple of the homology class
$[C_1]\in H_1(\Sg;\bZ)$ represented by the simple closed curve
$C_1$ with a fixed orientation depicted in Figure~\ref{fig:torelli}.
However, according to Johnson~\cite{Johnson80}, we have
$$
\tau(\varph)=(x_1\land y_1+\ldots +x_{g-1}\land y_{g-1})\land [C_1]
$$
where $x_1. y_1,\ldots ,x_{g-1}, y_{g-1}$ is a symplectic basis
of the homology group of the left hand subsurface of $\Sg$
obtained by cutting $\Sg$ along the simple closed curve $F_2$
depicted in Figure~\ref{fig:torelli}. It follows that
$$
C \tau(\varph)=2(g-1)[C_1] \ .
$$
This checks the equality~\eqref{eq:nutau}, for some constant $a$.

To determine this constant, it is enough to compute the
value of the cohomology class $[(\tilde\varph\inv)^*\nu-\nu]$
on the homology class represented by the oriented simple
closed curve $E$ depicted in Figure~\ref{fig:torelli}.
Observe that the homology class
$\tilde\varph\inv_*[E]-[E]\in H_1(\Sigma_g^{00};\bZ)$
can be represented by the oriented simple closed curve
$F_1$ also depicted in Figure~\ref{fig:torelli}.
Let $\widetilde D$ denote the right hand compact subsurface
of $\Sg$ obtained by cutting along $F_1$.
Thus $\widetilde D$ is diffeomorphic to a disk which contains
the original embedded disk $D$ in its interior.
Also let $\Sigma =\Sg\setminus\textrm{Int}\tilde D$.

Now we compute
\begin{align*}
[(\tilde\varph\inv)_*\nu-\nu]([E])&=\nu((\tilde\varph\inv)_*[E])-\nu([E])\\
&=\nu([F_1])
=\int_{F_1} \la-\la'\\
&=\int_{\partial \Sigma}\la-\int_{-\partial \widetilde D}\la'\\
&=\int_{\Sigma} d\la +\int_{\widetilde D}d\la'\\
&=\int_{\Sg} -\om =2-2g \ .
\end{align*}
Thus we can conclude that $a=-1$, completing the proof.
\end{proof}

\section{The second homology group of
$\Symp^{c,\delta}\Sigma_g^0$}\label{s:secondopen}

In this section, we generalize the claim of Theorem~\ref{th:hev}
about the second homology, i.~e.~for $k=1$, to the case of the
open surface $\Sigma_g^0$. We have the cohomology class
$$
j^*\tilde \al\in H^2(\Symp^{c,\delta}\Sigma_g^0;S^2_\bQ\bR)
$$
induced by the inclusion $j\colon\Sigma_g^0\ra\Sg$.
\begin{theorem}\label{th:h20}
The characteristic classes $e_1$ and $j^*\tilde \al$ induce a surjective
homomorphism
$$
H_2(\Symp^{c,\delta}\Sigma_g^0;\bZ)\lra \bZ\oplus S^2_\bQ\bR
$$
for any $g\geq 3$. For $g=2$, the class $j^*\tilde \al$ induces
a surjection
$$
H_2(\Symp^{c,\delta}\Sigma_2^0;\bZ)\lra S^2_\bQ\bR \ .
$$
\end{theorem}
The proof of this theorem, which occupies the rest of this section,
consists of a rather long and delicate argument. We first describe the
reason why the easy proof of Theorem~\ref{th:hev} in the case $k=1$,
which treated the case of
closed surfaces, does not work for Theorem~\ref{th:h20}, thereby
making the difficulty in the case of open surfaces explicit.

Consider the top extension
\begin{equation}\label{eq:scm}
1 \lra \Symp^{c}_{0}\Sigma_g^{0} \overset{i^c}{\lra}
\Symp^c\Sigma_g^{0} \overset{p}{\lra} \M1g \lra 1
\end{equation}
in the commutative diagram~\eqref{eqn:GG}. We have the following proposition,
contrasting with our discussion of the closed case in Section~\ref{s:main}.
\begin{proposition}\label{prop:fl2c}
For any $g\geq 2$ the image of the homomorphism
$$
(j^*\Flux)_*\colon H_2(\Symp_0^{c,\delta}\Sigma_g^0;\bZ)\lra
H_2(H^1(\Sg;\bR)^\delta;\bZ)\cong \La^2_{\bZ}H^1(\Sg;\bR)
$$
induced by the restriction $j^*\Flux$ of the flux homomorphism
to the subgroup $\Symp_0^{c}\Sigma_g^0$
is equal to the kernel of the natural intersection pairing
$$
\La^2_{\bZ}H^1(\Sg;\bR)\lra \bR \ .
$$
\end{proposition}
\begin{proof}
Recall that we have an extension
\begin{equation}\label{eq:hcss}
1\lra\Ham^c\Sigma_g^0\lra\Symp_0^c\Sigma_g^0\ \overset{j^*\Flux}{\lra}\
H^1(\Sg;\bR)\lra 1 \ ,
\end{equation}
where $\Ham^c\Sigma_g^0$ denotes the subgroup consisting of
Hamiltonian symplectomorphisms with compact supports
(see~\cite{MS} as well as~\cite{KM03}).
Also recall that there is a surjective homomorphism
$$
{\rm Cal}\colon \Ham^c \Sigma_g^0 \lra\bR
$$
called the (second) Calabi homomorphism (see~\cite{Calabi}).
Banyaga~\cite{Banyaga} proved that the kernel of this
homomorphism is perfect. Hence we have an isomorphism
$$
H_1(\Ham^{c,\delta}\Sigma_g^0;\bZ)\cong \bR \ .
$$
Now we consider the Hochschild--Serre exact sequence
\begin{align*}
H_2(\Symp_0^{c,\delta}\Sigma_g^0;\bZ)&\ \overset{j^*\Flux_*}{\lra}\
H_2(H^1(\Sg;\bR)^\delta;\bZ)\overset{\partial}{\lra}
H_1(\Ham^{c,\delta}\Sigma_g^0;\bZ)_{H^1_{\bR}}\cong\bR\\
&\lra H_1(\Symp_0^{c,\delta}\Sigma_g^0;\bZ)\
\overset{j^*\Flux_*}{\lra}\
H_1(H^1(\Sg;\bR)^\delta;\bZ)\lra 0
\end{align*}
of the group extension~\eqref{eq:hcss}. We proved in~\cite{KM03}
(Proposition~11 and Corollary~12) that the last homomorphism
$j^*\Flux_*$ in the above sequence is an isomorphism, and that
the boundary operator $\partial$ coincides with the intersection
pairing. The result follows.
\end{proof}

\begin{corollary}\label{cor:trivial}
The restriction of $\al\in H^2(\Symp^\delta\Sg;\bR)$ to the subgroup
$\Symp_0^{c,\delta}\Sigma_g^0$ is trivial.
\end{corollary}
\begin{proof}
This follows from Proposition~\ref{prop:fl2c} and the definition of
the cohomology class $\al$.
\end{proof}
Thus, in order to prove the non-triviality of $\al$ on the group
$\Symp^c\Sigma_g^0$, we must combine the roles of the two groups $\M1g$
and $\Symp_0^{c,\delta}\Sigma_g^0$. This contrasts sharply
with the case of closed surfaces treated in Section~\ref{s:main}.

Let $\{E_{p,q}^r\}$ be the Hochschild--Serre spectral sequence for the
integral homology of the extension~\eqref{eq:scm}.
This gives rise to two short exact sequences
\begin{equation}\label{eq:ss}
\begin{split}
0\lra \Ker& \lra H_2(\Symp^{c,\delta}\Sigma_g^0;\bZ)\lra
E^\infty_{2,0} \lra 0\\
& 0\lra E^\infty_{0,2}\lra\Ker \lra E^\infty_{1,1}\lra 0 \ ,
\end{split}
\end{equation}
where $E^\infty_{2,0} \subset H_2(\M1g;\bZ)$ concerns the first
Mumford-Morita-Miller class already discussed in~\cite{KM03}.
Proposition~\ref{prop:fl2c} shows that the image of the map
$$
E^\infty_{0,2}=\im\left(H_2(\Symp_0^{c,\delta}\Sigma_g^0;\bZ)
\ra H_2(\Symp^{c,\delta}\Sigma_g^0;\bZ)\right)\lra S^2_\bQ\bR
$$
defined by the Kronecker product with $\tilde\al$ is
precisely the kernel of the natural map $S^2_\bQ\bR\ra\bR$.
It remains to determine the $E^\infty_{1,1}$-term.
We have
\begin{equation}\label{eq:E1}
E^2_{1,1}=H_1(\M1g;H_1(\Symp_0^{c,\delta}\Sigma_g^0;\bZ))
\cong H_1(\M1g;H_1(\Sg;\bR))
\end{equation}
because, as was already mentioned above, it was proved in~\cite{KM03}
that $j^*\Flux$ induces an isomorphism
$H_1(\Symp_0^{c,\delta}\Sigma_g^0;\bZ)\cong H^1(\Sg;\bR)$.

We now recall a result which was essentially proved in~\cite{Morita89}.
Without repeating everything done in~\cite{Morita89}, we  want to give
a precise statement and proof of what is needed in the sequel.
We use the Lickorish generators for the mapping class group $\M1g$,
denoted $\la_i, \mu_i, \nu_i$ as in Figure~1 of~\cite{Morita89}, and a
symplectic basis $x_1,\cdots,x_g,y_1,\cdots,y_g$ of $H_1(\Sg;\bZ)$.
In particular
\begin{equation}\label{eq:action}
\nu_i(x_i)=x_i-y_i+y_{i+1},\quad \mu_i(x_i)=x_i-y_i
\end{equation}
which we record here for later use.
\begin{proposition}\label{prop:h1m1}
For any $g\geq 2$, we have an isomorphim
$$
H_1(\M1g;H_1(\Sg;\bZ))\cong \bZ \ .
$$
Furthermore, for any $i=1,2,\cdots,g-1$, the element
$$
c=\nu_i\otimes x_i-\mu_i\otimes x_i+\mu_{i+1}\otimes x_{i+1}
$$
is a $1$-cycle of $\M1g$ with twisted coefficients in $H_1(\Sg;\bZ)$
and it represents a generator of the above infinite cyclic group.
\end{proposition}
\begin{proof}
The group extension
$$
1\lra \pi_1\Sg\lra \Msg\lra \Mg\lra 1 \ ,
$$
where $\Msg$ denotes the mapping class group of $\Sg$ {\it relative}
to a basis point, yields the Hochschild--Serre exact sequence
\begin{equation}\label{eq:msg}
H_2(\Mg;H)\lra (H\otimes H)_{\Mg}\lra H_1(\Msg;H)
\lra H_1(\Mg;H)\lra 0 \ .
\end{equation}
Here and henceforth $H$ is a shorthand for $H_1(\Sg;\bZ)$.
It is easy to see that the intersection pairing induces an isomorphism
$(H\otimes H)_{\Mg}\cong\bZ$ and the element $x_1\otimes y_1$, for
example, represents a generator. It was proved in~\cite{Morita89}
that
$$
H^1(\Msg;H^1(\Sg;\bZ))\cong Z,\quad
H_1(\Mg;H^1(\Sg;\bZ))\cong \bZ/(2g-2)\bZ \ .
$$
If we apply the crossed homomorphism $f\colon\Msg\ra H^1(\Sg;\bZ)$
given in the above cited paper, which detects a generator of the
above infinite cyclic group, to the element $x_1\otimes y_1$
considered as a cycle of $\Msg$ with coefficients in $H$, we obtain
$$
f(x_1)(y_1)=2-2g \ .
$$
On the other hand, $c$ is a cycle because
$$
\partial c=x_i-\nu_i(x_i)-x_i+\mu_i(x_i)+x_{i+1}-\mu_{i+1}(x_{i+1})
=0
$$
by~\eqref{eq:action}, and it was shown that
\begin{equation}\label{eq:k}
f(c)=f(\nu_i)(x_i)-f(\mu_i)(x_i)+f(\mu_{i+1})(x_{i+1})=1 \ ,
\end{equation}
see~\cite{Morita89} for details. In view of the exact sequence~\eqref{eq:msg},
we can conclude that $H_1(\Msg;H)\cong \bZ$. Finally, it is easy to deduce
from the central group extension $0\ra \bZ\ra \M1g\ra\Msg\ra 1$ that we have an
isomorphism $H_1(\M1g;H)\cong H_1(\Msg;H)$. This finishes the proof.
\end{proof}

Going back to the $E^2_{1,1}$-term in~\eqref{eq:E1}, we have
$$
E^2_{1,1}\cong H_1(\M1g;H\otimes\bR)\cong\bR
$$
by Proposition~\ref{prop:h1m1}. Now consider the exact sequence
$$
E^2_{3,0}\cong H_3(\M1g;\bZ)\ \overset{d^2}{\lra} \
E^2_{1,1}\cong \bR\lra E^3_{1,1}=E^\infty_{1,1}\lra 0 \ .
$$
Harer~\cite{Harer92} determined the third stable
rational cohomology group
$$
\lim_{g\to\infty} H^3(\M1g;\bQ)
$$
of the mapping class group to be trivial\footnote{There is now a final
result on the stable cohomology of $\Mg$ due to Madsen and Weiss~\cite{MW}.}.
It follows that $E^2_{3,0}$ is a finite group for all sufficiently large $g$.
Hence we can conclude that $E^\infty_{1,1}\cong \bR$ for such $g$.
It is natural to expect that this $\bR$ will recover the missing $\bR$ in
$E^\infty_{0,2}$ so that we obtain the surjectivity of the $\tilde\al$-factor
in Theorem~\ref{th:h20}. It turns out that this is indeed the case, and below
we shall give a proof of this fact which does not use Harer's result mentioned
above. Before doing so, we have to prepare some general facts concerning
group (co)homology in small degrees.

Consider a group extension
\begin{equation}\label{eq:ext}
1\lra K\lra G\lra Q\lra 1 \ ,
\end{equation}
and
suppose we are given a $1$-cycle $c=\sum_i q_i\otimes u_i\in Z_1(Q;H_1(K))$
of the group $Q$ with coefficients in the abelianization $H_1(K)$ of $K$
considered as a natural $Q$-module, where $q_i\in Q$ and $u_i\in H_1(K)$.
\begin{lemma}\label{lem:cc}
For any choices of lifts $\tilde q_i\in G$ of the $q_i$
and representatives $k_i\in K$ with $[k_i]=u_i$, the element
$$
\tilde c=\sum_{i} \left\{(\tilde q_i, k_i)+(\tilde q_i k_i,\tilde q_i^{-1})
-(\tilde q_i,\tilde q_i^{-1})-(\id,\id)\right\}+ d
$$
is a $2$-cycle of $G$, where $d$ is a $2$-chain of the group $K$ such that
$$
\partial d=\sum_i\left\{(\tilde q_i k_i\tilde q_i^{-1})-(k_i)\right\}
\ .
$$
Furthermore, $p_*([\tilde c])=0\in H_2(Q;\bZ)$, where $p\colon G\ra Q$ denotes
the projection, and in the short exact sequence
\begin{align*}
0\lra &E^\infty_{0,2}\lra
\Ker\left(H_2(G;\bZ)\overset{p_*}{\lra} H_2(Q;\bZ)\right)\\
& \lra E^\infty_{1,1}\
\left(\cong H_1(Q;H_1(K))/d^2(E^2_{3,0})\right)\lra 0
\end{align*}
arising from the Hochschild--Serre spectral sequence of~\eqref{eq:ext} the
class $[\tilde c]\in \Ker$ is a lift of $[c]\in H_1(Q;H_1(K))$.
\end{lemma}
\begin{proof}
Since $c$ is a cycle by the assumption, we have
$$
\partial c=\sum_i (q_i(u_i)-u_i)=0 \in H_1(K) \ .
$$
It follows that there exists a $2$-chain $d\in C_2(K;\bZ)$
with the property described in the statement of the lemma.
Then a direct computation shows that $\partial \tilde c=0$.
Clearly $p_*(\tilde c)=0$ and it is easy to check
the rest of the required assertions.
\end{proof}

The following can be proved by a standard argument in the cohomology theory
of groups, see~\cite{Brown}.
\begin{lemma}\label{lem:ch}
Let $G$ be a group and $M$ a $G$-module. Assume we have a
$G$-invariant skew-symmetric bilinear pairing
$$
\iota\colon M\times M\lra A \ ,
$$
where $A$ is an abelian group with trivial $G$-action, and
we are given two crossed homomorphisms
$$
f_i\colon G\lra M\quad (i=1,2)
$$
so that $f_i(gh)=f_i(g)+ g_* f_i(h)$ for all $g,h\in G$.
Then the assignment
$$
G\times G\ni (g,h)\longmapsto\iota(f_1(g),g_*f_2(h))\in A \ ,
$$
which we denote by $f_1\cdot f_2$, is a $2$-cocycle of $G$ with values
in $A$ and its cohomology class in $H^2(G;A)$ depends only on the cohomology
classes $[f_i]\in H^1(G;M)$ of the crossed homomorphisms $f_i$.
Furthermore, $f_2\cdot f_1$ is cohomologous to $f_1\cdot f_2$,
so that $[f_2\cdot f_1]=[f_1\cdot f_2]\in H^2(G;A)$.
\end{lemma}

Now in the situation of Lemma~\ref{lem:ch}, we consider the case where
$G=\Symp^c\Sigma_g^0$, $M=H^1(\Sg;\bR)$ and $\iota$ is the intersection
pairing. Then we have three crossed homomorphisms
$$
j^*\widetilde{\Flux},\ \widetilde{\Flux}_c,\ p^*k_\bR:
\Symp^c\Sigma_g^0\lra H^1(\Sg;\bR)
$$
and, by the definition of $\al$, we have
$j^*\al=[j^*\widetilde{\Flux}\cdot j^*\widetilde{\Flux}]$.
\begin{proposition}\label{prop:fluxc}
In the above notation, we have
$[\widetilde{\Flux}_c\cdot\widetilde{\Flux}_c]=0$
and
$$
j^*\al=2 [p^*k_\bR\cdot\widetilde{\Flux}_c]-p^* e_1 \ .
$$
\end{proposition}
\begin{proof}
Define a map
$$
\widetilde{\Cal}\colon\Symp^c\Sigma_g^0\lra \bR
$$
by the formula
$$
\widetilde{\Cal}(\varph)=\int_{\Sigma_g^0}(\varph\inv)^*\la\land \la
\ .
$$
The restriction of this map to the subgroup $\Ham^c\Sigma_g^0$ is a
homomorphism, which, suitably normalized, is called the (second) Calabi
homomorphism (see~\cite{Calabi}). In our previous paper~\cite{KM03}, we
examined how the restriction of $\widetilde{\Cal}$ to the subgroup
$\Symp^c_0\Sigma_g^0$ fails to be a homomorphism. Here we extend this
discussion to the whole group $\Symp^c\Sigma_g^0$. For any two elements
$\varph, \psi\in \Symp^c\Sigma_g^0$, we claim that
\begin{equation}\label{eq:cal}
\widetilde{\Cal}(\varph\psi)
= \widetilde{\Cal}(\varph)+\widetilde{\Cal}(\psi)+
\widetilde{\Flux}_c(\varph)\cdot (\varph\inv)^*\widetilde{\Flux}_c(\psi)
\ .
\end{equation}
This is because
\begin{align*}
\widetilde{\Cal}(\varph\psi)=&\int_{\Sigma_g^0}
((\varph\psi)\inv)^*\la\land\la\\
=& \int_{\Sigma_g^0} (\varph\inv)^*((\psi\inv)^*\la-\la)\land\la
+(\varph\inv)^*\la\land\la\\
=& \int_{\Sigma_g^0}(\varph\inv)^*\widetilde{\Flux}_c(\psi)\land
((\varph\inv)^*\la-\widetilde{\Flux}_c(\varph))+\widetilde{\Cal}(\varph)\\
=& \widetilde{\Flux}_c(\varph)\cdot (\varph\inv)^*\widetilde{\Flux}_c(\psi)
+\int_{\Sigma_g^0} ((\psi\inv)^*\la-\la)\land\la +\widetilde{\Cal}(\varph)\\
=& \widetilde{\Flux}_c(\varph)\cdot (\varph\inv)^*\widetilde{\Flux}_c(\psi)
+\widetilde{\Cal}(\psi)+\widetilde{\Cal}(\varph) \ .
\end{align*}
Equation~\eqref{eq:cal} implies that the $2$-cocycle
$\widetilde{\Flux}_c\cdot \widetilde{\Flux}_c$ of the group $\Symp^c\Sigma_g^0$
is a coboundary. Hence $[\widetilde{\Flux}_c\cdot
\widetilde{\Flux}_c]=0$ as claimed.

The definition of $\al$ implies
$$
j^*\al= [j^*\widetilde{\Flux}\cdot j^*\widetilde{\Flux}]
$$
so that, by Theorem~\ref{th:ffk}, we can write
\begin{equation}\label{eq:all}
j^*\al=
[\widetilde{\Flux}_c\cdot\widetilde{\Flux}_c]+[\widetilde{\Flux}_c\cdot
p^*k_\bR]+[p^*k_\bR\cdot \widetilde{\Flux}_c]+p^*[k_\bR\cdot k_\bR] \ .
\end{equation}
By Lemma~\ref{lem:ch} we have
$[\widetilde{\Flux}_c\cdot p^*k_\bR]=[p^*k_\bR\cdot
\widetilde{\Flux}_c]$,
and it was proved in~\cite{Morita89b} that $[k_\bR\cdot k_\bR]=-e_1$.
If we substitute these relations in \eqref{eq:all}, we obtain the desired
identity.
\end{proof}

\begin{proof}[Proof of Theorem~\ref{th:h20}]
In view of Proposition~\ref{prop:fl2c} and the discussion following
Corollary~\ref{cor:trivial}, it suffices to show that there exist $2$-cycles
of the group $\Symp^c\Sigma_g^0$ such that the evaluations of $j^*\al$ on
them have as values any real number. To construct such $2$-cycles, we use
the $1$-cycle $c$ of $\M1g$ with coefficients in $H$ described in
Proposition~\ref{prop:h1m1}, which represents a generator of 
$H_1(\M1g;H)\cong\bZ$.
In order to adapt to the situation here, we consider the dual cycle
$$
c^*=\nu_i\otimes y^*_i-\mu_i\otimes y^*_i+\mu_{i+1}\otimes y^*_{i+1}
$$
with coefficients in $H^1(\Sg;\bZ)$, where the
$x^*_1,\cdots,x^*_g,y^*_1,\cdots,y^*_g$ denote the dual basis of 
$H^1(\Sg;\bZ)$.
If $k\colon\M1g\ra H^1(\Sg;\bZ)$ is a crossed homomorphism which represents
a generator of $H^1(\M1g;H^1(\Sg;\bZ))$, then we have
\begin{equation}\label{eq:kc}
k(c^*)=k(\nu_i)\cdot y^*_i-k(\mu_i)\cdot y^*_i+k(\mu_{i+1})\cdot
y^*_{i+1}=1 \ .
\end{equation}
For any real number $r\in\bR$, we consider the $1$-cycle
$$
c^*_r=\nu_i\otimes r y^*_i-\mu_i\otimes r y^*_i+\mu_{i+1}\otimes r y^*_{i+1}
$$
of $\M1g$ with coefficients in $H^1(\Sg;\bR)$.
Now we apply Lemma~\ref{lem:cc} in the case where
$G=\Symp^c\Sigma_g^0$, $K=\Symp^c_0\Sigma_g^0$, and $Q=\M1g$.
We proved in~\cite{KM03} that the flux homomorphism induces an isomorphism
$$
\Flux\colon H_1(\Symp^c_0\Sigma_g^0;\bZ)\cong H^1(\Sg;\bR).
$$
Therefore,
$c^*_r$ can be considered as a $1$-cycle of $\M1g$ with coefficients
in the abelianization of $\Symp^c_0\Sigma_g^0$.
Hence, if we choose elements
$$
\tilde{\nu}_i, \tilde{\mu}_i\in \Symp^c\Sigma_g^0,
\quad \varph^r_i\in \Symp^c_0\Sigma_g^0
$$
such that
\begin{equation}\label{eq:kf}
p(\tilde{\nu}_i)=\nu_i,\quad p(\tilde{\mu}_i)=\mu_i,
\quad \Flux(\varph^r_i)= \Flux_c(\varph^r_i)=r y^*_i \ ,
\end{equation}
where $p\colon\Symp^c\Sigma_g^0\ra\M1g$ denotes the projection as before,
then
\begin{align*}
\tilde{c}^*_r=& (\tilde{\nu}_i, \varph^r_i)+
(\tilde{\nu}_i\varph^r_i,\tilde{\nu}^{-1}_i)-
(\tilde{\nu}_i,\tilde{\nu}^{-1}_i)-(\id,\id)+\\
& (\tilde{\mu}_i, \varph^r_i)+
(\tilde{\mu}_i\varph^r_i,\tilde{\mu}^{-1}_i)-
(\tilde{\mu}_i,\tilde{\mu}^{-1}_i)-(\id,\id)+\\
&  (\tilde{\mu}_{i+1}, \varph^r_{i+1})+
(\tilde{\mu}_{i+1}\varph^r_{i+1},\tilde{\mu}^{-1}_{i+1})-
(\tilde{\mu}_{i+1},\tilde{\mu}^{-1}_{i+1})-(\id,\id)+d
\end{align*}
is a $2$-cycle of $\Symp^c\Sigma_g^0$, where $d$ is a $2$-chain of the
group $\Symp^c_0\Sigma_g^0$ such that
$$
\partial d=(\tilde{\nu}_i\varph^r_i\tilde{\nu}^{-1}_i)-(\varph^r_i)
+(\tilde{\mu}_i\varph^r_i\tilde{\mu}^{-1}_i)-(\varph^r_i)
+(\tilde{\mu}_{i+1}\varph^r_{i+1}\tilde{\mu}^{-1}_{i+1})-(\varph^r_{i+1}) \ .
$$
Now we claim that
\begin{equation}\label{eq:cr}
j^*\al(\tilde{c}^*_r)=2r \ ,
\end{equation}
which will finish the proof of the theorem.
To show this, observe first that $p^*e_1(\tilde{c}^*_r)=0$ because
clearly $p_*(\tilde{c}^*_r)=0$. Hence, by Proposition~\ref{prop:fluxc},
$$
j^*\al(\tilde{c}^*_r)=2 [p^*k\cdot\widetilde{\Flux}_c](\tilde{c}^*_r)
\ .
$$
Observe that
$$
[p^*k\cdot\widetilde{\Flux}_c]\left((\tilde{\nu}_i\varph^r_i,
\tilde{\nu}^{-1}_i)-(\tilde{\nu}_i,\tilde{\nu}^{-1}_i)\right)=0 \ ,
$$
because $p^*k(\tilde{\nu}_i\varph^r_i)=p^*k(\tilde{\nu}_i)$ and
$\varph^r_i$ acts trivially on the homology of $\Sigma_g^0$.
The same is true for two other similar terms. Since $d$ is a $2$-chain
of $\Symp^c_0\Sigma_g^0$, the evaluation of $[p^*k\cdot\widetilde{\Flux}_c]$
on it vanishes. Keeping in mind equations~\eqref{eq:kf} and~\eqref{eq:kc},
we can now conclude that
$$
j^*\al(\tilde{c}^*_r)=2 [p^*k\cdot\widetilde{\Flux}_c](\tilde{c}^*_r)
=2r p^*k(\tilde{c}^*)=2r \ .
$$
This proves~\eqref{eq:cr} and hence the theorem.
\end{proof}

\section{Proof of the main results}\label{s:proofs}

In this section we give the proofs of the main results described
in Section~\ref{s:main}.

\begin{proof}[Proof of Theorem~\ref{th:hev}]
Here we follow the argument of~\cite{Morita2} and of our previous
paper~\cite{KM03} to prove the non-triviality of the cup products
$e_1^k\tilde\al^{(\ell)}$. For this we first observe that, similar
to the class $e_1$, the class $\tilde\al$ is stable, with respect to
$g$, and also that it is primitive in the following sense. For each
$k$, consider the genus $kg$ surface
$\Sigma_{kg,1}=\Sigma_{kg}\setminus \Int D^2$ with one boundary component
as the boundary connected sum
$$
\Sigma_{kg,1}=\Sigma_{g,1}\ \natural \cdots\natural\ \Sigma_{g,1}
$$
of $k$ copies of $\Sigma_{g,1}=\Sg\setminus \Int D^2$.
This induces a homomorphism
\begin{equation}\label{eqn:fk}
f_k\colon \Symp^c \Sigma_g^0\times\cdots\times \Symp^c \Sigma_g^0
\lra \Symp^c \Sigma_{kg}^0
\end{equation}
from the direct product of $k$ copies of the group $\Symp^c \Sigma_g^0$ to
$\Symp^c \Sigma_{kg}^0$.
Under this homomorphism we have the equality
$$
f_{k}^*(\tilde\al)=\tilde\al\times 1\times\cdots\times 1+\cdots
+1\times\cdots\times 1\times \tilde\al \ ,
$$
which follows easily from the definition of $\tilde\al$. Now we can combine
Theorem~\ref{th:h20} with the above property to show the required assertion
in the theorem. This finishes the proof.
\end{proof}

\begin{proof}[Proof of Theorem~\ref{th:om2}]
The fact that the ideal generated by $\om_0\land H_1(T^{2g};\bR)$
is contained in the kernel of $\Flux^*$ has already been proved in
Proposition~\ref{prop:v}. To show that $\Ker \Flux^*$ is precisely
this ideal, we use the decomposition~\eqref{eq:irrep} of $H^k(T^{2g};\bR)$
into irreducible summands given in Section~\ref{s:main}. It is easy to
see that the quotient of this module divided by the ideal generated by
$\om_0\land H_1(T^{2g};\bR)\subset H_3(T^{2g};\bR)$ is precisely
$$
\bR\oplus \bigoplus_{k=1}^{g} [1^k]_\bR \ .
$$
The $\bR$-summand in degree $2$ corresponds to the class $\alpha$
and its non-triviality has already been shown in Theorem~\ref{th:hev}.
Hence, to prove the assertion, it remains to show that $\Flux^*([1^k])$
is non-trivial for any $k\leq g$. The highest weight vector of the
irreducible representation $[1^k]$ is $x_1\land\cdots\land x_k$, where
$x_1,\cdots,x_g,y_1,\cdots,y_g$ is a symplectic basis of $H_1(\Sg;\bZ)$
as before. Now the definition of the flux homomorphism implies that, for
any $i$, there exists an element $\varph_i\in\Symp_0\Sg$ such that
$\Flux(\varph_i)=x_i^*$. Here we can choose the support of $\varph_i$ to
be contained in an arbitrarily small neighbourhood of a simple closed curve
which represents the homology class $x_i$. Then the $k$ elements
$\varph_1,\cdots,\varph_k$ mutually commute because their supports are
disjoint. Hence they form a cycle of $\Symp_0\Sg$ supported on a
$k$-dimensional torus, and the cohomology class
$\Flux^*(x_1\land\cdots\land x_k)\in H^{k}(\Symp_0^\delta\Sg;\bR)$
takes a non-zero value (namely $1$) on this cycle. This completes the proof.
\end{proof}

\begin{proof}[Proof of Theorem~\ref{th:ham}]
We consider the extension
$$
1\lra \Ham\Sg\lra \Symp_0\Sg\overset{\Flux}{\lra}H^1(\Sg;\bR)\lra 1 \ .
$$
Let $\{E^r_{p,q}\}$ be the Hochschild--Serre spectral sequence for its
homology. Since $\Ham\Sg$ is perfect by Thurston~\cite{Thurston0}, see
also~\cite{Banyaga}, we have
$$
E^2_{1,1}\cong H_1(H^1(\Sg;\bR)^\delta;H_1(\Ham^\delta\Sg;\bZ))=0 \ .
$$
Hence the differential $d^2\colon E^2_{3,0}\ra E^2_{1,1}$ vanishes, so that
$$
E^3_{3,0}\cong E^2_{3,0}\cong H_3(H^1(\Sg;\bR)^\delta;\bZ)
\cong \La_\bZ^3 H^1(\Sg;\bR) \ .
$$
Similarly $E^2_{2,1}=0$, so that the differential $d^2\colon E^2_{2,1}\ra 
E^2_{0,2}$
vanishes and
$$
E^3_{0,2}\cong E^2_{0,2}\cong
H_0(H^1(\Sg;\bR)^\delta;H_2(\Ham^\delta\Sg;\bZ))
\cong H_2(\Ham^\delta\Sg;\bZ)_{H^1_\bR} \ .
$$
However, $E^\infty_{3,0}=E^4_{3,0}$ is equal to the image of the homomorphism
$$
\Flux_*\colon H_3(\Symp_0^\delta\Sg;\bZ)\lra
H_3(H^1(\Sg;\bR)^\delta;\bZ) \ .
$$
Now we can conclude that the exact sequence
$$
E^\infty_{3,0}=E^4_{3,0}\lra E^3_{3,0}\lra E^3_{0,2}=E^2_{0,2}
$$
yields an exact sequence
\begin{equation}
H_3(\Symp_0^\delta\Sg;\bZ)\lra H_3(H^1(\Sg;\bR)^\delta;\bZ)\lra
H_2(\Ham^\delta\Sg;\bZ)_{H^1_\bR} \ .
\label{eq:eh}
\end{equation}
By the same computation as above in the dual context of cohomology,
we obtain an exact sequence
\begin{equation}
H^2(\Ham^\delta\Sg;\bR)^{H^1_\bR}\lra H^3(H^1(\Sg:\bR)^\delta;\bR)
\overset{\Flux^*}{\lra} H^3(\Symp_0^\delta\Sg;\bR) \ .
\label{eq:ech}
\end{equation}
Proposition~\ref{prop:v} (see also Theorem~\ref{th:om2})
implies that the continuous cohomology classes
$$
\tilde\om_0\land H_1(\Sg;\bR)\subset H^3(H^1(\Sg:\bR)^\delta;\bR)
$$
vanish under the homomorphism $\Flux^*$ so that they can be
lifted to elements of $H^2(\Ham^\delta\Sg;\bR)^{H^1_\bR}$
by~\eqref{eq:ech}. Now consider the cycles
$$
\om_0 \land H^1(\Sg;\bR)\subset \La^3_\bZ H^1(\Sg;\bR)\cong
H_3(H^1(\Sg;\bR)^\delta;\bZ) \ ,
$$
where $\om_0\in\La^2_\bZ H^1(\Sg;\bR)$, and also their images
in $H_2(\Ham^\delta\Sg;\bZ)_{H^1_\bR}$ in the exact sequence~\eqref{eq:eh}.
If we consider the Kronecker products of these cycles with the above lifted
cohomology classes, we can conclude that $\om_0 \land H^1(\Sg;\bR)$ maps
injectively into $H_2(\Ham^\delta\Sg;\bZ)_{H^1_\bR}$. This completes the proof.
\end{proof}

\begin{remark}
It would be interesting to obtain explicit group cocycles of $\Ham\Sg$ which
represent the above degree two cohomology classes. There should also be a
relation to the work of Ismagilov~\cite{Is}. We shall pursue this elsewhere.
\end{remark}

\section{Further discussion}\label{s:misc}

As in Section~\ref{s:trans}, let $e, v\in H^2(\ES^{\delta}\Sg;\bR)$ be
the Euler class and the transverse symplectic class, respectively.
By analogy with the definition $e_1=\pi_*(e^2)$, where
$$
\pi_*\colon H^4(\ES^{\delta}\Sg;\bR)\lra H^2(\Symp^\delta\Sg;\bR)
$$
denotes the integration over the fiber, we can define a cohomology class
$$
v_1 \in H^2(\Symp^\delta\Sg;\bR)
$$
by setting $v_1=\pi_*(ev)$. After we conjectured that $v_1$ is a linear
combination of the two classes $\al$ and $e_1$, Kawazumi~\cite{Kawazumi}
kindly provided a proof. More precisely, he pointed out that the
contraction formula, Theorem~6.2 of~\cite{KM01}, can be adapted to the
case of the cohomology class $\al$ of the group $\ES^{\delta}\Sg$, and
that the following equality holds:
\begin{equation}
\al=-\pi_*\left((e+v)^2\right)=-e_1-2 v_1 \ .
\label{v1}
\end{equation}


Since we know by~\cite{KM03} that $e^2\not=0$, we could also apply
integration over the fiber to the cohomology class $e^2v$ in order to
obtain some more cohomology. However, $e^2 v$ vanishes in
$H^6(\ES^\delta\Sg;\bR)$ by the Bott vanishing theorem.



A more promising approach to find more cohomology for the
symplectomorphism groups is the following.
Consider the extension
\begin{equation}
1\lra\Symp_0\Sg\lra\Symp\Sg\overset{p}{\lra}\Mg\lra 1 \ .
\label{eq:ssm}
\end{equation}
On the one hand, Theorem~\ref{th:om2} shows that we have an injection
\begin{equation}
\bigoplus_{k=1}^g [1^k]_\bR \subset H^*(\Symp_0^\delta\Sg;\bR) \ .
\label{eq:1k}
\end{equation}
On the other hand, Looijenga~\cite{Looijenga} determined the stable
cohomology $H^*(\Mg;V)$ of the mapping class group with coefficients
in any irreducible representation $V$ of the algebraic group $Sp(2g;\bQ)$.
In particular, the cohomology groups $H^*(\Mg;[1^k])$ are highly non-trivial.
In the spectral sequence for the cohomology of the extension~\eqref{eq:ssm},
there are many non-trivial $E_2$-terms
$$
E^{p,q}_2=H^p(\Mg;H^q(\Symp_0^\delta\Sg;\bR)) \ .
$$
For example, it was proved in~\cite{Morita93} that 
$H^1(\Mg;[1^3]_\bR)\cong\bR$,
and an explicit computation using Looijenga's formula shows
$H^2(\Mg;[1^2]_\bR)\cong\bR$. Hence we have injections
\begin{align*}
&\bR\subset E^{1,3}_2=H^1(\Mg;H^3(\Symp_0^\delta\Sg;\bR)) \\
&\bR\subset E^{2,2}_2=H^2(\Mg;H^2(\Symp_0^\delta\Sg;\bR))
\end{align*}
for all sufficiently large $g$. It seems likely that these
two copies of $\bR$ survive to the $E_\infty$ term, so that they
define certain cohomology classes in $H^4(\Symp^\delta\Sg;\bR)$.
More generally, the summands~\eqref{eq:1k} and the non-trivial
cohomology groups $H^*(\Mg;[1^k]_\bR)$ should give rise to
infinitely many cohomology classes in $H^*(\Symp^\delta\Sg;\bR)$.
\begin{problem}
Prove that these cohomology classes are non-trivial.
\end{problem}

Next, there are completely different candidates for possible classes
in $H^*(\Symp^\delta\Sg;\bR)$ coming from the cohomology of the Lie
algebra $V_n$ of formal Hamiltonian vector fields on $\bR^{2n}$
first studied by Gelfand, Kalinin and Fuks in~\cite{GKF}. This Lie
algebra $V_n$ contains ${\mathfrak sp}(2n,\bR)$ as a subalgebra
consisting of vector fields corresponding to linear symplectomorphisms.
Let $\BG_2^\om$ be the Haefliger classifying space for the pseudogroup
of local symplectomorphisms of $\bR^2$ with respect to the standard
symplectic form. Then there is a natural homomorphism
\begin{equation}
H^*_c(V_1;{\frak sp}(2,\bR))\lra H^*(\BG_2^\om;\bR)
\label{eq:vbg}
\end{equation}
from the continuous cohomology of $V_1$ relative to the subalgebra
${\frak sp}(2,\bR)$ to the real cohomology group of $\BG_2^\om$.

There is also an obvious continuous mapping
$$
K(\ES^\delta\Sg,1)\lra \BG_2^\om
$$
which classifies the transversely symplectic codimension $2$
foliation on the classifying space for the group $\ES^\delta\Sg$,
that is the total space of the universal foliated $\Sg$-bundle over
$\BS^\delta\Sg$. This induces homomorphisms
\begin{equation}
H^*(\BG_2^\om;\bR)\lra H^*(\ES^\delta\Sg;\bR)
\overset{\pi_*}{\lra} H^{*-2}(\Symp^\delta\Sg;\bR) \ ,
\label{eq:bgs}
\end{equation}
where the last homomorphism is the integration along the fibre.
Combining~\eqref{eq:vbg} and~\eqref{eq:bgs} we obtain a homomorphism
\begin{equation}
H^*_c(V_1;{\frak sp}(2,\bR))\lra H^{*-2}(\Symp^\delta\Sg;\bR) \ .
\label{eq:}
\end{equation}

Now, Gelfand, Kalinin and Fuks~\cite{GKF} found a new cohomology class
in $H^7_c(V_1;{\frak sp}(2,\bR))$.
Later, Metoki~\cite{Metoki} extended their computation and found another
exotic class in $H^9_c(V_1;{\frak sp}(2,\bR))$. It seems to be widely
believed that there should exist infinitely many such exotic classes.
\begin{problem}
Study the cohomology classes in $H^*(\Symp^\delta\Sg;\bR)$ induced from
exotic classes in $H^{*+2}_c(V_1;{\frak sp}(2,\bR))$. In particular, prove
that the two elements in $H^5(\Symp^\delta\Sg;\bR)$ and 
$H^7(\Symp^\delta\Sg;\bR)$
induced from the classes found by Gelfand, Kalinin and Fuks and by Metoki,
are non-trivial.
\end{problem}

Recall that Harer~\cite{Harer85} proved that the homology groups of the mapping
class groups $\Mg$ stabilize with respect to the genus $g$ in a certain stable
range. In view of the fact that all the characteristic classes we introduced in
this paper are stable with respect to $g$, we would like to propose the 
following
problem, although it appears to be beyond the range of available techniques 
at the
moment:
\begin{problem}
Determine whether the homology groups of $\Symp^\delta \Sg$ stabilize
with respect to $g$, or not. In particular, is it true that
$$
H_2(\Symp^\delta\Sg;\bZ)\cong \bZ\oplus S^2_\bQ\bR
$$
for all $g\geq 3$ ?
\end{problem}

\section*{Appendix: Proof of Proposition~\ref{prop:la2}}

To prove Proposition~\ref{prop:la2}, observe first that the second exterior
power $\La^2_{\bZ}H^1(\Sg;\bR)$ over $\bZ$ is naturally isomorphic to the same
power $\La^2_{\bQ}H^1(\Sg;\bR)$ over $\bQ$ because $H^1(\Sg;\bR)$ is a uniquely
divisible abelian group. Choose a Hamel basis $a_\la\ (\la\in A)$ of $\bR$ as a
vector space over $\bQ$. Then we can write
$$
H^1(\Sg;\bR)=\sum_{\la} a_\la H^1(\Sg;\bQ) \ .
$$
Hence
$$
\La^2_{\bQ} H^1(\Sg;\bR)=\sum_{\la} a_\la \La^2_{\bQ}H^1(\Sg;\bQ)\oplus
\sum_{\la<\mu} a_\la H^1(\Sg;\bQ)\otimes a_\mu H^1(\Sg;\bQ) \ ,
$$
where we choose a total order in the index set $A$. Clearly, this is a
decomposition of $\Mg$-modules. It is easy to see that the intersection
pairing gives rise to an isomorphism
$$
\left(\La^2_{\bQ}H^1(\Sg;\bQ)\right)_{\Mg}\cong \bQ \ .
$$
We also have
$$
\left(H^1(\Sg;\bQ)\otimes H^1(\Sg;\bQ)\right)_{\Mg}
\cong
\left(S^2H^1(\Sg;\bQ)\oplus \La^2H^1(\Sg;\bQ)\right)_{\Mg}\\
\cong \bQ \ .
$$
Here we have used the fact that
$$
\left(S^2H^1(\Sg;\bQ)\right)_{\Mg}=0 \ ,
$$
which is true because the action of $\Mg$ on $S^2H^1(\Sg;\bQ)$ factors
through that of the algebraic group $Sp(2g,\bQ)$, and $S^2H^1(\Sg;\bQ)$
is a non-trivial irreducible $Sp(2g,\bQ)$-module.

Thus we obtain an isomorphism
$$
\left(\La^2_{\bZ} H^1(\Sg;\bR)\right)_{\Mg}\cong
\left(\sum_{\la} a_\la\otimes a_\la \bQ\right)\oplus
\left(\sum_{\la<\mu} a_\la\otimes a_\mu\bQ\right) \ .
$$
It is easy to see that the right-hand side of the above expression
can be naturally identified with $S^2_\bQ\bR$. This completes the proof.

\bibliographystyle{amsplain}

\end{document}